\documentclass[reqno]{amsart}

\usepackage{amsmath,amssymb,latexsym}
\usepackage[english]{babel}

\usepackage{color}

\usepackage{setspace}

\newcommand{\field}[1]{\ensuremath{{\mathbb #1}}}
\newcommand{\Quat}{\field{H}}
\newcommand{\R}{\field{R}}
\newcommand{\C}{\field{C}}
\newcommand{\N}{\field{N}}

\newcommand{\eqskip}{\quad}

\newcommand{\Harm}{\ensuremath{\mathrm{Harm}}}

\newcommand{\dm}{\,d}

\DeclareMathOperator{\sign}{sign}

\newtheorem{theorem}{Theorem}[section] 
\newtheorem{lemma}[theorem]{Lemma}
\newtheorem{prop}[theorem]{Proposition}
\newtheorem{cor}[theorem]{Corollary}

\theoremstyle{remark}
\newtheorem{remark}[theorem]{Remark}

\newcommand{\ltwoab}{\ensuremath{L_2^{(\alpha,\beta)}}}

\numberwithin{equation}{section}

\begin{document}

\author{Yu. I. Lyubich}
\title{ Lower bounds for projective designs, cubature formulas 
and related isometric embeddings}
 
\begin{abstract}
  Yudin's lower bound \cite{yudin} for the spherical designs is
  generalized to the cubature formulas on the projective spaces 
  over a field $\field{K}\subset \{\R,\C,\Quat\}$ and thus
  to isometric embeddings $l_{2;\field{K}}^m\rightarrow l_{p;\field{K}}^n$ 
  with $p\in 2\N$. For large $p$ and in some other situations 
  this is essentially better than those known before.
  

\end{abstract}

\maketitle

\section{Introduction}
\label{sec:intro}

In the theory of spherical and projective designs some important lower
bounds were obtained
\cite{bannaihoggarsymmetric,delsartegoethalsseidel} by maximization of
the functional
\[
  D(f) = \frac{f(1)}{c_0[f]},\eqskip f\in K_l,
\]
where $K_l$ is the set of nonnegative on $(-1,1)$ 
nonzero polynomials $f$, deg$f\leq l$, and
\begin{equation}
  \label{eq:1}
  c_0[f] = \int_{-1}^1 f(t)\omega_{\alpha,\beta}(t)\dm t,\eqskip
  \omega_{\alpha,\beta}(t) = (1-t)^\alpha(1+t)^\beta.
\end{equation}
Here the numbers $l\in\N$ and $\alpha$, $\beta > -1$ depend on the design. 

Obviously, $\sup\left\{ D(f) : f\in K_l\right\} = \sup\left\{f(1) :
  f\in K_l, c_0[f]=1\right\}$. The solution to the latter linear
programming problem is classical, the extremal polynomial $f_{\rm max}$ is
unique and can be expressed in terms of the Jacobi polynomials, see
\cite{szegogabor}, Section 7.7.1. For the designs of cardinality $n$
this yields
\begin{equation}
  \label{eq:2}
  n\geq \tau_{\alpha,\beta}f_{\rm max}(1) =
  \tau_{\alpha,\beta}\max\left\{D(f): f\in K_l\right\},
\end{equation}
where
\begin{equation}
  \label{eq:3}
  \tau_{\alpha,\beta} = \int_{-1}^1 \omega_{\alpha,\beta}(t)\dm t.
\end{equation}

We denote by $\ltwoab(-1,1)$ the space of complex-valued
measurable functions $f$ on $(-1,1)$ such that
\[
  \|f\|^2 \equiv \int_{-1}^1 |f(t)|^2
  \omega_{\alpha,\beta}(t)\dm t<\infty .
\]
The corresponding Jacobi polynomials $P_k(t)$ constitute an orthogonal
basis in $\ltwoab$, so that
\begin{equation}
  \label{eq:4}
  f(t) = \sum_{k=0}^\infty \nu_k c_k[f] P_k(t),\eqskip f\in \ltwoab,
\end{equation}
where
\begin{equation}
  \label{eq:5}
  c_k[f] = \int_{-1}^1 f(t) P_k(t)\omega_{\alpha,\beta}(t)\dm t,\eqskip
  \nu_k=1/\|P_k\|^2.
\end{equation}
 \eqref{eq:4} The Jacobi-Fourier series \eqref{eq:4} converges to $f$ in
$\ltwoab(-1,1)$. The coefficient $c_0[f]$ in \eqref{eq:5} coincides
with that of \eqref{eq:1} since $P_0(t)\equiv 1$, according to the
usual standardization
\[
  \deg P_k = k,\eqskip P_k(1) = \binom{\alpha+k}{k} \eqskip
  (k=0,1,2,\ldots).
\]
For the same reason $\nu_0 = 1/\tau_{\alpha,\beta}$.

The linear programming bound \eqref{eq:2} can be extended 
to the set $K_{l,l'}$, $l'>l$, of the polynomials $f\neq 0$, $\deg
f\leq l'$, such that $f(t) \geq 0$ for $|t|\leq 1$ and $c_k[f]\leq 0$
for $l+1 \leq k \leq l'$. In this way Boyvalenkov and Nikova 
\cite{boyv95,boyvnikov94,boyvnikov97} obtained a series of new 
concrete lower bounds for the projective designs. For the spherical designs
Yudin \cite{yudin} considered the limit case $l' = \infty$. Its class 
$K_{l,\infty}$ consists of all nonnegative nonzero continuous functions
$f(t)$, $|t|\leq 1$, such that $c_k[f]\leq 0 $ for all $k\geq l+1$.  
A suitable choice of a function $f\in K_{l,\infty}$ yields 
a lower bound asymptotically better than classical one 
that comes from \eqref{eq:2}.

In the present paper we generalize Yudin's result on the projective
designs and even the cubature formulas on the projective spaces 
$\field{K}P^{m-1}$, where $\field{K}\subset \{\R,\C,\Quat\}$. (Recall that 
$\field{H}$ is the standard notation for the quaternion field.) The extension 
of the linear programming bound from the projective designs 
to the general cubature formulas is technically simple 
but important since the latter are equivalent to the isometric embeddings 
$l_{2;\field{K}}^m\rightarrow l_{p;\field{K}}^n$, $p\in 2\N$.
(See \cite{lyushatquaternions} and references therein.) 
Note that with the standard inner product $(x,y)$ the space
$l_{2;\field{K}}^m$ is Euclidean, its unit sphere is $S = S^{\delta m-1}$,
 where $\delta = \delta(\field{K})$ such that
$\delta (\R)=1,\delta (\C) = 2, \delta (\Quat) = 4$ .

From now on we assume $m\geq2$, $p\in 2\N$, and denote by $\Phi_\field{K}(m,p)$
the space of complex-valued functions $\phi(x)$, $x\in S$, satisfying
the following conditions, see \cite{lyushatquaternions, lyushatpolyfunc}.
\begin{enumerate}
\item $\phi=\psi|S$, where $\psi$ is a homogeneous polynomial of
  degree $p$ (``$p$-forms'') on the space $\field{K}^m\equiv \R^{\delta
    m}$;
\item $\psi$ is invariant in the sense that
  \[
     \psi(w\alpha)  = \psi(w)\eqskip (w\in \field{K}^m, \alpha\in
     \field{K},\eqskip |\alpha|=1).
  \]
\end{enumerate}
A fortiori, $\phi(x\alpha)=\phi(x)$ that allows us to naturally
transfer $\phi$ to the projective space $\field{K}P^{m-1}$. However,
we will consider $\phi$ on $S$ which is equivalent but more
elementary. In this setting a projective {\em cubature formula of
  index $p$ on $S$} is
\begin{equation}
  \label{eq:6}
  \int \phi\dm \sigma = \sum_{i=1}^n \phi(x_i)\rho_i,\eqskip \phi\in\Phi_\field{K}(m,p),
\end{equation}
where $\sigma$ is the normalized Lebesgue measure on $S$, the nodes
$x_i\in S$ are projectively distinct and the weights $\rho_i$ are
positive. (Note that $\sum\rho_i = 1$ automatically by the restriction
of $(x,x)^{p/2}$ to $S$.) In the case of equal $\rho_i$ the set
$\left\{x_i\right\}_1^n$ is nothing but a projective $p/2$-design,
c.f. \cite{hoggardesignsproj}.

\section{Basic theory}
\label{sec:basic-theory}

First of all, we have the decomposition
\begin{equation}
  \label{eq:2-1}
  \Phi_\field{K}(m,p) = \sum_{k=0}^{p/2} \Harm_\field{K}(m,2k)
\end{equation}
where the space $\Harm_\field{K}(m,2k)$ consists of restrictions
to $S$ of the invariant harmonic $2k$-forms. Regarding to the inner
product
\[
  (\psi_1,\psi_2) = \int\bar\psi_1\psi_2\dm\sigma
\]
the decomposition \eqref{eq:2-1} is orthogonal.

For any orthonormal basis $\{\phi_{ki}\}_{i=1}^{d_{m,2k}}$ of
$\Harm_\field{K}(m,2k)$ the {\em addition formula}
\begin{equation}
  \label{eq:2-2}
  \sum_{s=1}^{d_{m,2k}} \overline{\phi_{ks}(x)}\phi_{ks}(y) =
    b_{m,k}P_k(xy) \eqskip (x,y\in S)
\end{equation}
holds with
\begin{equation} 
 \label{eq:2-3}
  b_{m,k} = \tau_{\alpha,\beta}\nu_k P_k(1), \eqskip \alpha = 
\frac{\delta(m-1)-2}{2},\eqskip \beta=\frac{\delta-2}{2},
\end{equation}
and
\begin{equation}
  \label{eq:2-3p}
  xy = 2|(x,y)|^2 - 1,
\end{equation}
see \cite{hoggarpreprint,koornwinder,lyushatquaternions,muller}. Later
on we operate only with $\alpha,\beta$ given by \eqref{eq:2-3}. 

Now let $X$ be a finite nonempty subset of $S$, and let $A(X)$ be its angle set,i.e.
\[
  A(X) = \{xy : x,y\in X, x\neq y \}.
\]
The addition formula easily implies the following
\begin{lemma}
  \label{lem:2-1}
  Let the series
  \[
     \sum_{k=0}^\infty a_k P_k(t)
  \]
  converge to a function $f(t)$ for every $t\in A(X)$ and for
  $t=1$. Then
  \begin{equation}
    \label{eq:2-4}
    \sum_{x,y\in X} f(xy)\bar\lambda(x)\lambda(y) = \sum_{k=0}^\infty
    a_k b_{m,k}^{-1}\sum_{s=1}^{d_{m,2k}} \left\lvert \sum_{x\in X}
      \phi_{ks}(x)\lambda(x) \right\rvert^2,
  \end{equation}
  where $\lambda$ is an arbitrary function $X\rightarrow\C$.
\end{lemma}

With $\lambda(x)\equiv 1$ formula \eqref{eq:2-4} plays a fundamental
role in the design theory
\cite{delsartegoethalsseidel,hoggardelsartespaces,levenshtein}. In the
context of cubature formulas we need \eqref{eq:2-4} with arbitrary
$\lambda(x)>0$, $\sum\lambda(x)=1$, c.f.\cite{lyushatquaternions}, \S5. 
Also note that, in contrast to those which are quoted above, now we have
to apply \eqref{eq:2-4} to the non-polynomial functions $f\in
K_{l,\infty}$. It is possible because of
\begin{lemma}
  \label{lem:2-2}
  The Jacobi-Fourier series of any function $f\in K_{l,\infty}$
  converges to $f(t)$ for all $t\in[-1,1]$.
\end{lemma}
\begin{proof}
  Since $f(t)$ is continuous, its Jacobi-Fourier series at $t=1$ is
  summable to $f(1)$ by a Cesaro method, see \cite{szegogabor},
  Theorem 9.13. Therefore, it is summable to $f(1)$ by the Abel
  method, see \cite{hardy}, Theorem 5.5. Hence, this series converges to
  $f(1)$ since $c_k{f}\leq 0$ for $k\geq l+1$. It remains to refer to
  Theorem 7.32.1 from \cite{szegogabor} which states that
  \begin{equation}
    \label{eq:2-4p}
    \max_{|t|\leq 1} \left\lvert P_k(t)\right\rvert = P_k(1)
  \end{equation}
  if $\max(\alpha,\beta)\geq -1/2$. The latter is fulfilled
  because of \eqref{eq:2-3} and $m\geq 2$.
\end{proof}

\begin{cor}
  \label{cor:2.2p}
  Formula \eqref{eq:2-4} is true for every $f\in K_{l,\infty}$ with
  $a_k =\nu_k c_k[f]$, $k\geq 0$.
\end{cor}

\begin{remark}
  \label{rem:2-3}
  In \cite{yudin} the absolute convergence of the corresponding series is
  mentioned without proof. The proof of Lemma \ref{lem:2-2} shows that 
in our situation the convergence is absolute and uniform.
\end{remark}

Now we can prove the following linear programming bound.

\begin{prop}
  \label{prop:2-4}
  The inequality
  \begin{equation}
    \label{eq:2-5}
    n\geq \tau_{\alpha,\beta} \sup\left\{ D(f) : f\in
      K_{p/2,\infty}\right\}
  \end{equation}
  holds for any projective cubature formula of shape \eqref{eq:6}.
\end{prop}
\begin{proof}
  We have
  \begin{equation}
    \label{eq:2-6}
    \sum_{i=1}^n \phi(x_i)e_i = 0,\eqskip \phi\in
    \Harm_\field{K}(m,2k),\eqskip 1\leq k\leq p/2.
  \end{equation}
  Applying Corollary \ref{cor:2.2p} and Lemma \ref{lem:2-1} to $f\in
  K_{p/2,\infty}$, $X=\{x_i\}_1^n$ and $\lambda(x)=\rho(x)$, $x\in X$,
  we obtain
  \[
    f(1) \sum_{x\in X}\rho^2(x) \leq \sum_{x,y\in X}
    f(xy)\rho(x)\rho(y) \leq a_0b_{m,0}^{-1}
    \left(\sum\rho(x)\right)^2 .
  \]
  Indeed, on the left side of \eqref{eq:2-4} all summands are
  $\geq0$. On the right side the summands are $\leq 0$ for $k\geq p/2 + 1$ 
  and vanish for $1\leq k \leq p/2$ by
  \eqref{eq:2-6}. It remains to recall that $\sum\rho(x) = 1$,
  therefore, $\sum\rho^2(x)\geq n^{-1}$; on the other hand, $a_0
  b_{m,0}^{-1} = c_0[f]/\tau_{\alpha,\beta}$  since $b_{m,0} =1$.
\end{proof}

\begin{remark}
  \label{rem:2-4p}
   The inequality \eqref{eq:2-5} implies
  \[
    n \geq \tau_{\alpha,\beta} \sup\{D(f) : f\in K_{p/2,l'}\}, \eqskip l' \geq p/2, 
  \]
  since $K_{p/2,l'}\subset K_{p/2,\infty}$. (For $l' =l$ we set $K_{l',l} = K_l$.)
\end{remark}

Since with any given $m,p$ a projective cubature formula exists (or,
equivalently, there exists an isometric embedding $l_{2;
  \field{K}}^m\rightarrow l_{p; \field{K}}^n$ ), we have
\begin{cor}
  \label{cor:2-5}
   $\sup\{D(f) : f\in K_{p/2,\infty}\} < \infty$.
\end{cor}

The supremum in question is unknown but a ``good'' test function can
be constructed using the ``convolution''
\begin{equation}
  \label{eq:2-7}
  \int g(xu)h(uy)\dm\sigma(u)\eqskip (x,y\in S)
\end{equation}
of two suitable functions $g(t)$ and $h(t)$, $-1\leq t \leq 1$, 
c.f. \cite{yudin}. 
\begin{lemma}
  \label{lem:2-6}
  For any $e\in L_1^{(\alpha,\beta)}(-1,1)$ the function $u\mapsto
  e(xu)$, $u\in S$, belongs to $L_1(S,\sigma)$ for every $x\in S$ and
  \begin{equation}
    \label{eq:2-8}
    \int e(xu)\dm\sigma(u) = \frac{1}{\tau_{\alpha,\beta}} \int_{-1}^1
    e(t)\omega_{\alpha,\beta}(t)\dm t
  \end{equation}
\end{lemma}
\begin{proof}
This follows by calculation in spherical coordinates.
\end{proof}
\begin{cor}
  \label{cor:2-7}
  With $g,h\in L_2^{(\alpha,\beta)}(-1,1)$ the integral \eqref{eq:2-7}
  exists for all $x,y\in S$.
\end{cor}

Since any ordered pair $x',y'\in S$ with $x'y' = xy$ can be obtained
from $x,y$ by an isometry of $l_{2;\field{K}}^m$, the integral
\eqref{eq:2-7} depends on $xy$ only. Thus, we have a function
$(g*h)(t)$, $-1\leq t\leq 1$, such that
\begin{equation}
  \label{eq:2-9}
  (g*h)(xy) = \int g(xu)h(uy)\dm\sigma(u).
\end{equation}
In particular, for $x=y$ \eqref{eq:2-9} yields
\begin{equation}
  \label{eq:2-10}
  (g*h)(1) = \int g(xu)h(xu)\dm\sigma(u) =
  \frac{1}{\tau_{\alpha,\beta}} \int_{-1}^1
  g(t)h(t)\omega_{\alpha,\beta}(t)\dm t,
\end{equation}
by \eqref{eq:2-8}. Moreover, applying the Schwartz inequality to
\eqref{eq:2-9} and using \eqref{eq:2-8} again we obtain
\begin{equation}
  \label{eq:2-11}
  \sup_t |(g*h)(t)|\leq \frac{1}{\tau_{\alpha,\beta}} \|g\|\cdot\|h\|.
\end{equation}
By this inequality and bilinearity, the convolution $g*h$ determines a
continuous mapping $(L_2^{(\alpha,\beta)})^2 \rightarrow L_\infty$. 

\begin{lemma}
  \label{lem:2-8}
  With $g,h\in L_2^{(\alpha,\beta)}$ the function $(g*h)(t)$ is
  continuous, and the series
  \[
    (g*h)(t) = \sum_{k=0}^\infty (\nu_k^2/b_{m,k})c_k[g] c_k[h]P_k(t)
  \]
converges uniformly.
\end{lemma}

\begin{proof}
  Let
  \[
    g_N = \sum_{j=0}^N \nu_j c_j[g]P_j, \eqskip h_N(t) = \sum_{k=0}^N
   \nu_k c_k[h]P_k .
  \]
  Then
  \[
    g_N * h_N = \sum_{k=0}^N(\nu_k^2/b_{m,k}) c_k[g]c_k[h]P_k
  \]
  since
  \[
    P_j * P_k = b_{m,k}^{-1} P_k\delta_{jk}
  \]
  by the addition formula. Since $g_N\rightarrow g$ and
  $h_N\rightarrow h$ $(N\rightarrow\infty)$ in
  $L_2^{(\alpha,\beta)}$, we obtain $g_N *h_N\rightarrow g*h$ uniformly.
  Thus, the limit function is continuous.
\end{proof}

\begin{cor}
  \label{cor:2-9}
  $c_k[g*h] = \nu_k c_k[g]c_k[h]/b_{m,k} = c_k[g] 
c_k[h]/\tau_{\alpha,\beta}P_k(1)$.
\end{cor}

\section{A function $f_l\in K_{l,\infty}$}
\label{sec:3}

Recall that all roots of every $P_k$ are simple and lie on $(-1,1)$.
The roots of the derivative $P'_k$ alternate them, so they are also simple and 
lie on $(-1,1)$. Now we introduce a function $f_l$ by setting 
\begin{eqnarray}
  \label{eq:3-1}
  f_l = g*h, & g(t)=
  \begin{cases}
    P_r(t) - P_r(\xi), & t\geq\xi \\
    0, & t<\xi 
  \end{cases}, &
  h(t) =
  \begin{cases}
    1, & t\geq \xi, \\
    0, & t<\xi
  \end{cases}
\end{eqnarray} 
where $r=l+1$ and $\xi$ is the largest root of $P'_r$. We have to 
verify that $f_l\in K_{l,\infty}$.

By Lemma \ref{lem:2-8} $f_l$ is continuous. The inequality $f_l\geq0$
follows from \eqref{eq:2-9} since $h\geq 0$ and $g\geq 0$. The former
is obvious, the latter is true since $g(\xi)=0$, $g(1)\geq 0$ and
$g'(t)\neq 0$ for $\xi < t\leq 1$. Moreover, $f_l(1)>0$ by
\eqref{eq:2-10}, thus, $f_l\neq 0$. It remains to prove that
$c_k[f_l]\leq 0$ for $k\geq r$. In \cite{yudin} a rather
complicated vector analysis on $\R^m$ was used at this point. 
We manage without a
generalization of this technique to $\C^m$ and $\Quat^m$ by dealing
with the corresponding Jacobi polynomials.

Our starting point is the differential equation 
\begin{equation}
\label{eq:3-2}
\Delta_i\equiv (\omega_{\alpha+1,\beta+1}P_i')' +
i(i+\lambda)\omega_{\alpha,\beta}P_i =0,\eqskip i\geq 0,
\end{equation}
where $\lambda=\alpha+\beta+1$, see \cite{szegogabor},
formula(4.2.1). Note that $\lambda\geq 0$ by \eqref{eq:2-3}.
From \eqref{eq:3-2} it follows that
\begin{equation*}
  \begin{split}
    0 & = \int_\xi^1 (P_r\Delta_k - P_k\Delta_r)\dm t = \\
      & = (k-r)(k+r+\lambda) \int_\xi^1 \omega_{\alpha,\beta} P_r P_k\dm t + 
    \int_\xi^1\{\omega_{\alpha+1,\beta+1} (P_r P_k'-P_k P_r'\}'\dm t = \\
      & = (k-r)(k+r+\lambda)\int_\xi^1 \omega_{\alpha,\beta} P_r P_k\dm t -
    (\omega_{\alpha+1,\beta+1}P_r P_k')(\xi)
\end{split}
\end{equation*}
since $\omega_{\alpha+1,\beta+1}(1)=0$, $P_r'(\xi)=0$. For $k\neq r$
we obtain
\[
  \int_\xi^1 P_r P_k\omega_{\alpha,\beta}\dm t =
  \frac{(\omega_{\alpha+1,\beta+1}P_r P_k')(\xi)}{(k-r)(k+r+\lambda)}.
\]
This formula extends to $r=0$ since $P_0(t)\equiv 1$, so
$P_0'(\xi)=0$. Thus,
\[
  \int_\xi^1 P_k \omega_{\alpha,\beta}\dm t = 
     \frac{(\omega_{\alpha+1,\beta+1}P_k')(\xi)}{k(k+\lambda)},
\]
and then
\[
  \int_\xi^1 P_r P_k \omega_{\alpha,\beta}\dm t = 
    \frac{k(k+r)P_l(\xi)}{(k-r)(k+r+1)}\int_\xi^1
    P_k\omega_{\alpha,\beta}\dm t.
\]
As a result,
\[
  c_k[g] = \int_\xi^1 g P_k \omega_{\alpha,\beta}\dm t =
  \frac{r(r+\lambda)P_r(\xi)}{(k-r)(k+r+\lambda)} \int_\xi^1
  P_k\omega_{\alpha,\beta}\dm t =
  \frac{r(r+\lambda)P_r(\xi)}{(k-r)(k+r+\lambda)}c_k[h],
\]

and, by Corollary \ref{cor:2-9},
\begin{equation}
  \label{eq:3-3}
  c_k[f_l] =
  \frac{r(r+\lambda)P_r(\xi)}{(k-r)(k+r+\lambda)P_k(1)\tau_{\alpha,\beta}}
    (c_k[h])^2 \eqskip (k\neq r).
\end{equation}
Since $P_k(1)>0$, formula \eqref{eq:3-3} yields 
$\sign c_k[f_l] = \sign(P_r(\xi)),\eqskip k>r$.
But $\sign P_r(\xi) =-1$ since $\xi$ lies in between two largest roots
of $P_r(t)$ and $P_r(1)>0$. Thus, $c_k[f_l]<0$ for $k>r$. In addition,
$c_r[f_l]=0$ since $c_r[h]=0$. The latter follows from \eqref{eq:3-2}
with $i=r$ by integration over $[\xi,1]$.

In conclusion we note that $\xi$ in \eqref{eq:3-1} is actually the
largest root of $P_l^{(\alpha+1,\beta+1)}(t)$, see \cite{szegogabor},
formula (4.21.7).

\section{Main Theorem}
\label{sec:4}

Now we are in position to prove the following

\begin{theorem}
  \label{thm:4-1}
  The number $n$ of nodes of every projective cubature formula of
  index $p$ on $S^{\delta m-1}$ satisfies the inequality
  \begin{equation}
    \label{eq:4-1}
    n\geq
    \frac{\Gamma(\alpha+2)\Gamma(\beta+1)}
{\Gamma(\alpha+\beta+2)F(-\beta,\alpha+1,\alpha+2,\varepsilon)} 
\left(\frac{1}{\varepsilon}\right)^{\delta(m-1)/2},
  \end{equation}
  where $F$ is the hypergeometric function, the numbers $\alpha$ and
  $\beta$ are given by \eqref{eq:2-3}, $\varepsilon=(1-\xi)/2$, $\xi$ is
  the largest root of the Jacobi polynomial $P^{(\alpha+1,\beta+1)}_{p/2}(t)$.
\end{theorem}
\begin{proof}
  Using $f_{p/2}(t)$ as a test function in \eqref{eq:2-5} we get
  \[
    n\geq \frac{\tau_{\alpha,\beta}f_{p/2}(1)}{c_0[f_{p/2}]}.
  \]
  By \eqref{eq:2-10} and \eqref{eq:3-1} we have
  \[\tau_{\alpha,\beta}f_{p/2}(1) = 
  \int_\xi^1 gh\omega_{\alpha,\beta}\dm t = 
  \int_\xi^1 g\omega_{\alpha,\beta}\dm t = c_0[g].
  \]
  On the other hand, $c_0[f_{p/2}] = c_0[g]c_0[h]/\tau_{\alpha,\beta}$
  by Corollary \ref{cor:2-9}. Hence,
  \begin{equation}
    \label{eq:4-2}
    n\geq \frac{\tau_{\alpha,\beta}}{c_0[h]} = \frac{\int_{-1}^1
      (1-t)^\alpha(1+t)^\beta\dm t}{\int_\xi^1 (1-t)^\alpha(1+t)^\beta
      \dm t}.
  \end{equation}
  Now we substitute $t=1-2s$ into the numerator and $t=1-2\varepsilon s$
  into thte denominator. This yields \eqref{eq:4-1} since
  \[
    F(-\beta,\alpha+1,\alpha+2,\varepsilon ) =(\alpha+1)\int_0^1 s^\alpha
    (1-\varepsilon s)^\beta\dm s,
  \]
  (c.f. \cite{abramstegun}, formula (15.3.1)) and
  \[
    \int_0^1 s^\alpha(1-s)^\beta\dm s =
    \frac{\Gamma(\alpha+1)\Gamma(\beta+1)}{\Gamma(\alpha+\beta+2)}.
  \]
  (Also note that $\alpha+1=(\delta m-\delta)/2$ by \eqref{eq:2-3}.)
\end{proof}

\begin{remark}
  By substitution $t=2s^2-1$ in both integrals \eqref{eq:4-2} we
  obtain
  \[
    n\geq \frac{\int_0^1(1-s^2)^\alpha s^{2\beta+1}\dm
      s}{\int_\eta^1(1-s^2)^\alpha s^{2\beta+1}\dm s}, \eqskip
    \eta=\sqrt{(1+\xi)/2}.
  \]
  In particular, for $\field{K} = \R$ we have $\alpha = (m-3)/2, 
  \beta = -1/2$, see \eqref{eq:2-3}. Hence, 
  \begin{equation}
    \label{eq:4-3}
    n\geq \frac{\int_0^1 (1-s^2)^{(m-3)/2}\dm s}
  {\int_\eta^1(1-s^2)^{(m-3)/2}\dm s} \eqskip (\field{K} =\R),  
  \end{equation}
  where $\eta$ is the largest root of the polynomial
  \begin{equation}
    \label{eq:4-4}
    P_{p/2}^{\left((m-1)/2,1/2\right)}(2s^2-1)=\mathrm{const}
    \cdot P_{p+1}^{\left ((m-1)/2,(m-1)/2\right)}(s)/s,
  \end{equation}
  or, equivalently, of the Gegenbauer polynomial
  $C_{p+1}^{m/2}(s)$ (see \cite{szegogabor}, formulas (4.1.5)
  and (4.7.1).). In the case of antipodal spherical $(p+1)$-design the
  lower bound \eqref{eq:4-3} turns into (3) of \cite{yudin} up to the
  additional factor $2$ in the latter. Note that the factor 2 is just the
  degree of the natural mapping $S^{m-1}\rightarrow\R P^{m-1}$.
\end{remark}
\begin{remark}
 By \eqref{eq:2-3} we have $\alpha =m-2, \beta = 0$ for $\field{K} = \C$,     
 and $\alpha =2m-3, \beta = 1$ for $\field{K} = \Quat$. Accordingly, 
\eqref{eq:4-1} reduces to
  \begin{equation}
    \label{eq:4-5}
    n\geq \left(\frac{1}{\varepsilon}\right)^{m-1} \eqskip (\field{K} = \C) ,
  \end{equation}
 and to
  \begin{equation}
    \label{eq:4-6}
    n\geq \frac{1}{(2m-1)-(2m-2)\varepsilon}\left(\frac{1}{\varepsilon}\right)^{2m-2}
  \eqskip  (\field{K}=\Quat).
  \end{equation}
  In the real case the hypergeometric function
  in \eqref{eq:4-1} is not a polynomial of $\epsilon$.
\end{remark}

Now we denote by $N_\field{K}(m,p)$ the minimal number $n$ of nodes in
the cubature formula \eqref{eq:6} or, equivalently, the minimal $n$
such that there is isometric embedding $l_{2;\field{K}}^m\rightarrow
l_{p;\field{K}}^n$. In this notation Theorem \ref{thm:4-1} states that
\begin{equation}
  \label{eq:4-7}
  N_{\field{K}}(m,p) \geq 
\frac{\Gamma(\alpha+2)\Gamma(\beta+1)}
{\Gamma(\alpha+\beta+2)F(-\beta,\alpha+1,\alpha+2,\varepsilon)}\left(\frac{1}{\varepsilon}
\right)^{\delta(m-1)/2}.
\end{equation}
We will compare this result to the linear programming bound
\eqref{eq:2} with $l = p/2$. An explicit form of the latter is
\begin{equation}
  \label{eq:4-8}
  N_{\field{K}}(m,p)\geq \Lambda_{\field{K}}(m,q),\eqskip q=p/2,
\end{equation}
where
\begin{equation}
  \label{eq:4-9}
  \Lambda_{\field{K}}(m,q) =
  \begin{cases}
    \binom{m+q-1}{m-1}, & (\field{K}=\R); \\
    \binom{m+[q/2]-1}{m-1}\binom{m+[(q+1)/2]-1}{m-1}, &
    (\field{K}=\C);\\
    \frac{1}{2m-1}\binom{2m+[q/2]-2}{2m-2}\binom{2m+[(q+1)/2]-1}{2m-2}, & 
   (\field{K}=\Quat),
  \end{cases}
\end{equation}
while
\begin{equation}
  \label{eq:4-10}
  N_\field{K}(m,p)\leq  \Lambda_\field{K}(m,p),
\end{equation}
(See \cite{lyushatquaternions} and the references therein.)

\section{Asymptotic analysis}
\label{sec:5asymptotic-analysis}

From \eqref{eq:4-8} and \eqref{eq:4-9} it follows that
\begin{equation}
  \label{eq:5-1}
  N_\field{K}(m,p)\gtrsim
  \frac{p^{\delta(m-1)}}{\lambda_\field{K}(m)},\eqskip p\rightarrow\infty,
\end{equation}
where
\begin{equation}
  \label{eq:5-2}
  \lambda_\field{K}(m)=
  \begin{cases}
    2^{m-1}(m-1)!, & \field{K}=\R; \\
    2^{4(m-1)}(m-1)!^2, & \field{K}=\C; \\
    2^{8(m-1)}(2m-1)!(2m-2)!, & \field{K}=\Quat.
  \end{cases}
\end{equation}
or, in an unified form,
\begin{equation}
  \label{eq:5-3}
  \lambda_\field{K}(m) = \frac{\Gamma(\delta m/2)\Gamma(\delta(m-1)/2+1)}{\Gamma(\delta /2} 
\cdot 2^{2\delta(m-1)} = 
\frac{\Gamma(\alpha+\beta+2)\Gamma(\alpha+2)}{\Gamma(\beta+1)} \cdot 2^{2\delta(m-1)}.
\end{equation}
As to \eqref{eq:4-7}, $\varepsilon$ is the only parameter depending on
$p$. (Of course, $\varepsilon$ also depends on $m$.) By definition,
$\varepsilon=(1-\xi)/2 = \sin^2(\theta/2)$ where
$\theta=\arccos\xi$. This $\theta$ is the smallest root of the
polynomial $P^{(\alpha+1,\beta+1)}_{p/2}(\cos\theta)$. By Theorem
8.1.2 from \cite{szegogabor} we have $\theta\sim 2j_{\alpha+1,1}/p$
where $j_{\alpha+1,1}$ is the smallest positive root of the Bessel's
function $J_{\alpha+1}(z)$. Therefore, $\varepsilon\sim
j_{\alpha+1,1}^2/p^2$, and \eqref{eq:4-7} yields 
\begin{equation}
  \label{eq:5-3p}
  N_\field{K}(m,p) \gtrsim
  \frac{\Gamma(\alpha+2)\Gamma(\beta+1)}{\Gamma(\alpha+\beta+2)} \cdot
  \frac{p^{\delta(m-1)}}{j_{\alpha+1,1}^{\delta(m-1)}}, \eqskip p\rightarrow\infty,
\end{equation}
since $\varepsilon\rightarrow 0$, $F(\cdot,\cdot,\cdot,0)=1$.\rm {\em This
estimate is better than \eqref{eq:5-1}} because of
\begin{prop}
  \label{prop:5-1}
  The inequality
  \begin{equation}
    \label{eq:5-4}
    j_{\alpha+1,1}^{\delta(m-1)}<
    \frac{\Gamma(\alpha+2)\Gamma(\beta+1)}{\Gamma(\alpha+\beta+2)} \lambda_\field{K}(m)
  \end{equation}
  holds for all $m\geq 2$, except for the case $m=2$, $\delta=1$, 
   when \eqref{eq:5-4} changes for an equality.
\end{prop}
\begin{proof}
  By \eqref{eq:5-3} the inequality \eqref{eq:5-4}is equivalent to
  \begin{equation}
    \label{eq:5-5}
    j_{\alpha+1,1}^{\delta(m-1)} < \Gamma(\alpha+2)^2\cdot 2^{2\delta(m-1)}
  \end{equation}
  We set $\alpha+1=\nu$, so that $\delta(m-1)=2\nu$, and
  \eqref{eq:5-5} takes the form
  \begin{equation}
    \label{eq:5-6}
    j_{\nu,1}^{2\nu} < \Gamma(\nu+1)^2\cdot 16^\nu.
  \end{equation}
  The number $\nu$ is positive integer or half-integer, $\nu\geq 1/2$, and $\nu = 1/2$
  if and only if $m=2$, $\delta=1$. In this case $j_{\nu,1}=\pi$ since 
$J_{1/2}(z)$ is
proportional to $\sin z/\sqrt{z}$. On the other hand, 
$\Gamma(3/2)^2\cdot 16^{1/2} = \pi$        
as well. Thus, \eqref{eq:5-6} changes for an equality.

Now let $\nu\geq 1$. By the inequality $j_{\nu,1} < \sqrt{2(\nu+1)(\nu+3)}$ 
(see \cite{watson}, Section 15.3) it suffices to prove that
  \begin{equation}
    \label{eq:5-7}
    (\nu+1)^\nu(\nu+3)^\nu \leq \Gamma(\nu+1)^2\cdot 8^\nu .
  \end{equation}
 By Stirling's lower bound the inequality \eqref{eq:5-7} follows from
  \[
    \left(1+\frac{1}{\nu}\right)^\nu\left(1+\frac{3}{\nu}\right)^\nu <
    2\pi\nu\left(\frac{8}{e^2}\right)^\nu 
  \]
  A fortiori, \eqref{eq:5-7} follows from
  \[
    2\pi\nu\left(\frac{8}{e^2}\right)^\nu > e^4.
  \]
  But the latter is indeed true if $\nu\geq\nu_0$ where $\nu_0$ is a
  unique root of the equation $2\pi\nu(8/e^2)^\nu = e^4$. It is easy
  to see that $\nu_0< 6$, so \eqref{eq:5-7} is valid for $\nu\geq
  6$. For $\nu < 6$, i.e. $\nu= 1,3/2,2,\ldots,5,11/2$,  
  the inequality \eqref{eq:5-7} can be checked numerically.
\end{proof}

The inequalities \eqref{eq:5-1} and \eqref{eq:5-3p} can be rewritten as
\begin{equation}
  \label{eq:5-8}
  \liminf_{p\rightarrow\infty} p^{-\delta(m-1)}N_\field{K}(m,p) \geq 1/\lambda_\field{K}(m)
\end{equation}
and
\begin{equation}
  \label{eq:5-9}
  \liminf_{p\rightarrow\infty} p^{-\delta(m-1)}N_\field{K}(m,p) 
\geq \frac{\Gamma(\alpha+2)\Gamma(\beta+1)}{\Gamma(\alpha+\beta+2)} 
j_{\alpha+1,1}^{-\delta(m-1)}
\end{equation}
respectively. By Proposition \ref{prop:5-1} the quotient
$\kappa_\field{K}(m)$ of the lower bounds \eqref{eq:5-8} and
\eqref{eq:5-9} is less than $1$, except for the case $m=2$,
$\delta=1$. Moreover, $\kappa_\field{K}(m)$ exponentially decays 
as $m\rightarrow\infty$, c.f. \cite{yudin} for the spherical designs. 

\begin{prop}
  \label{prop:5-2}
  Asymptotically,
  \begin{equation}
    \label{eq:5-10}
    \kappa_\field{K}(m) \sim \frac{1}{\pi\delta
      m}\left(\frac{e}{4}\right)^{\delta(m-1)}, \eqskip m\rightarrow\infty.
  \end{equation}
\end{prop}
\begin{proof}
  In the same notation as before we have
  \[
    \kappa_\field{K}(m) =
    \frac{j_{\nu,1}^{2\nu}}{\Gamma(\nu+1)^2\cdot 16^\nu}.
  \]
  Using the relation
  $j_{\nu,1} = \nu + O(\sqrt[3]{\nu}), \eqskip \nu\rightarrow\infty$ ,
  (see \cite{watson}, Section 15.83) and Stirling's asymptotic formula
   we obtain
  \[
    \kappa_\field{K}(m) \sim
    \frac{1}{2\pi\nu}\left(\frac{e}{4}\right)^{2\nu}
  \]
  that is equivalent to \eqref{eq:5-10} since $2\nu \sim \delta m$.
\end{proof}

\begin{remark}
From \eqref{eq:4-10} the asymptotic upper bound
\begin{equation}
  \label{eq:5-11}
  \limsup_{p\rightarrow\infty} p^{-\delta(m-1)}N_\field{K}(m,p)\leq
  2^{\delta(m-1)}/\lambda_\field{K}(m)
\end{equation}
follows. We see that there is an exponential gap between \eqref{eq:5-11} 
and \eqref{eq:5-9} as $m\rightarrow\infty$. 
Indeed,the quotient of these bounds is
\begin{equation}
  \label{eq:5-12}
  2^{\delta(m-1)}\kappa_\field{K}(m) \sim \frac{1}{\pi\delta
    m}\left(\frac{e}{2}\right)^{\delta(m-1)}, \eqskip m\rightarrow\infty.
\end{equation}
\end{remark}

\section{The case $m=2$}
\label{sec:6}

In this case  we discuss the real, complex and quaternion situation separately.

\subsection{$\field{K}=\R$}
\label{sec:6a}
{\em Then the inequalities \eqref{eq:4-7} and \eqref{eq:4-8} are both the
equalities,so they coincide.} Indeed, $N_\R(2,p) = p/2+1$, according to
\cite{lyuvas,reznick}, and, on the other hand, 
$\Lambda_\R(2,q) = q+1 = p/2+1$ by \eqref{eq:4-9}. Furthermore, 
in the real case \eqref{eq:4-7} is equivalent to \eqref{eq:4-3}.
For $m=2$ this yields $N_\R(2,p)\geq \pi /2\arccos\eta = p/2 + 1$.
Indeed, in this context $\eta$ is the largest root of the Gegenbauer polynomial
$C_{p+1}^1(s)$ = sin$(p+2)\theta$/sin$\theta$ where $\theta$ = arccos$s$.

\subsection{$\field{K}=\C$.}
\label{sec:6b}

By \eqref{eq:4-9}
\begin{equation}
  \label{eq:6-3}
  N_\C(2,p)\geq \left[\left(\frac{p}{4}+1\right)^2\right],
\end{equation}
c.f.\cite{konig}. On the other hand, our bound \eqref{eq:4-5} for $m=2$ is
\begin{equation}
  \label{eq:6-4}
  N_\C(2,p)\geq\left]\frac{2}{1-\xi_p}\right[ 
\end{equation}
where $\xi_p$ is the largest root of $P_{p/2}^{(1,1)}(t)$ and
$]\zeta[$ means the smallest ingeger $\geq\zeta$, $\zeta\in\R$. A
numerical evaluation shows that \eqref{eq:6-4} coincides with
\eqref{eq:6-3} for $p\leq 16$, but exceeds it for $18\leq p\leq 90$. Moreover, the 
difference $\Delta_\C(p)$ between the lower bounds
\eqref{eq:6-4} and \eqref{eq:6-3} is nondecreasing in this range, as
we see from the table
\begin{table}[h]
\begin{tabular}{|c|c|c|c|c|c|c|c|c|c|c|c|c|c|c|c|c|c|c|c|c|c|c|c|c|c|c|c|c|c|c|c|}
  \hline
  $p$ & $\leq 16$ & 18 -- 24 & 26,28 & 30,32 & 34, 36 & 38 & 40 & 42, 
  44 & 46 & 48 & 50 & 52 \\ \hline
  $\Delta_\C(p)$ & 0 & 1 & 2 & 3 &4 & 5 & 6 & 7 & 8 & 9 & 10 & 11 \\ \hline
\end{tabular}
\label{tab:16-52}
\end{table}

\begin{table}[h]
\begin{tabular}{|c|c|c|c|c|c|c|c|c|c|c|c|c|c|c|c|c|c|}
  \hline
   $p$ & 54 & 56 & 58 & 60 & 62 & 64 & 66 & 68 & 70
   & 72 & 74 & 76 & 78 & 80 & 82 & 84 \\ \hline
   $\Delta_\C(p)$ & 12 & 13 & 14 & 15 & 17 & 18 & 19 & 20 & 22 & 23 &
   25 & 26 & 28 & 29 & 31 & 32 \\ \hline
\end{tabular}
\end{table}

\begin{table}[h]
\begin{tabular}{|c|c|c|c|}
\hline
$p$ &  86 & 88 & 90 \\ \hline
$\Delta_\C(p)$ & 34 & 36 & 38 \\ \hline
\end{tabular}
\end{table}
The table  also shows that the ``derivative `` $\Delta'_\C(p) =
\Delta_\C(p)-\Delta_\C(p-2)$ is nondecreasing (rather slowly).

\subsection{\field{K}=\Quat}
\label{sec:6c}

We have
\begin{equation}
  \label{eq:6-5}
  N_\Quat(2,p) \geq \frac{1}{3}\binom{[p/2]+2}{2}\binom{[(p+2)/2]+3}{2}
\end{equation}
from \eqref{eq:4-8} and \eqref{eq:4-9}, but \eqref{eq:4-6} yields
\begin{equation}
  \label{eq:6-6}
  N_\Quat(2,p) \geq \left]\frac{4}{(2+\eta_p)(1-\eta_p)^2}\right[
\end{equation}
where $\eta_p$ is the largest root of $P_{p/2}^{(2,2)}(t)$. 

Comparing \eqref{eq:6-6} to \eqref{eq:6-5} one can see a small advantage of 
\eqref{eq:6-5} when $4\leq p \leq 20$, $p\neq 18$. Namely, for the difference 
$\Delta_\Quat(p)$ between
the lower bounds \eqref{eq:6-6} and \eqref{eq:6-5} we have 

\begin{table}[h]
\begin{tabular}{|c|c|c|c|c|c|c|c|c|c|c|}\hline
  $p$ & 2 & 4 & 8 & 8 & 10 & 12 & 14 & 16 & 18 & 20 \\ \hline
  $\Delta_\Quat(p)$ & 0 & -1 & -1 & -4 & -2 & -6 & -3 & -6 & 1 & -1 \\ \hline
\end{tabular}
\caption{}
\label{tab:6-7}
\end{table}
However, for $p \geq 22$ this difference increases rather rapidly:

\begin{table}[h]
\begin{tabular}{|c|c|c|c|c|c|c|c|c|c|c|c|c|c|c|c|}\hline
  $p$ & 22 & 24 & 26 & 28 & 30 & 32 & 34 & 36 & 38 & 40 & 42 & 44 & 46
  & 48 & 50 \\ \hline
  $\Delta_\Quat(p)$ & 12 & 14 & 35 & 42 & 75 & 90 & 138 & 165 & 231 &
  274 & 364 & 426 & 544 & 631 & 782 \\ \hline 
\end{tabular}
\caption{}
\label{tab:6-8}
\end{table}

Also, an interesting observable phenomenon is a regular oscillation of
$\Delta_\Quat'(p)$ in contrast to the monotonicity of
$\Delta_\C'(p)$. Indeed, in both tables \ref{tab:6-7} and \ref{tab:6-8} we have
\begin{equation}
  \label{eq:6-9}
  \sign\Delta_\Quat''(p) = (-1)^{p/2+1}.
\end{equation}
for the second difference. This can be conjectured for all $p$ 
as well as the results of observations above.

\bibliographystyle{plain}
\bibliography{lowbndprojdes}

\end{document}